\documentclass[12pt,twoside]{article}
\usepackage{times}
\usepackage{amsmath,amssymb}
\usepackage{amsthm}
\usepackage{mathrsfs}
\usepackage{color}
\usepackage{exscale}
\usepackage{relsize}
\usepackage{color}
\usepackage{fancyhdr}
\allowdisplaybreaks

\newcommand{\ef}{ \hfill $ \Box $ \vskip 3mm}
\newcommand{\be}{\begin{equation}}
\newcommand{\ee}{\end{equation}}
\newcommand{\bea}{\begin{eqnarray}}
\newcommand{\eea}{\end{eqnarray}}

\newcommand{\bR}{{\mathbb R}}
\newcommand{\bN}{{\mathbb N}}

\def\nn{\nonumber}

\def\ve{\varepsilon}

\def\q{\quad}
\def\qq{\qquad}
\def\th{\theta}

\def\Dl{\Delta}
\def\ve{\varepsilon}

\def\lt{\left}

\def\rt{\right}

\def\i{\infty}

\def\supp{\text{supp }}

\def\p{\partial}
\def\f{\frac}
\def\na{\nabla}
\def\al{\alpha}

\def\O{\Omega}

\def\s{\sqrt}
 \allowdisplaybreaks

\begin{document}
 \footskip=0pt
 \footnotesep=2pt
\let\oldsection\section
\renewcommand\section{\setcounter{equation}{0}\oldsection}
\renewcommand\thesection{\arabic{section}}
\renewcommand\theequation{\thesection.\arabic{equation}}
\newtheorem{theorem}{Theorem}[section]
\newtheorem{proposition}[theorem]{Proposition}
\newtheorem{remark}[theorem]{Remark}
\newtheorem{corollary}[theorem]{Corollary}

\title{Liouville theorem of axially symmetric Navier-Stokes equations with growing velocity at infinity}

\author{Xinghong Pan$^{a,}$\footnote{E-mail:xinghong{\_}87@nuaa.edu.cn}\ ,\quad Zijin Li $^{b,}$\footnote{E-mail:zijinli@smail.nju.edu.cn} \vspace{0.5cm}\\
\footnotesize $^a$Department of Mathematics, Nanjing University of Aeronautics and Astronautics, Nanjing 211106, China.\\
 \footnotesize $^b$College of Mathematics and Statistics, Nanjing University of Information Science $\&$ Technology,\\
  \footnotesize Nanjing 210044, China.\\
\vspace{0.5cm}
}

\date{}

\maketitle

\centerline {\bf Abstract} In the paper \cite{KNSS:1}, the authors make the following conjecture: {\it any bounded ancient mild solution of the 3D axially symmetric Navier-Stokes equations is constant.} And it is proved in the case that the solution is swirl free. Our purpose of this paper is to improve their result by allowing that the solution can grow with a power smaller than 1 with respect to the distance to the origin. Also, we will show that such a power is optimal to prove the Liouville type theorem since we can find counterexamples for the Navier-Stokes equations such that the Liouville theorem fails if the solution can grow linearly.
\vskip 0.3 true cm

\vskip 0.3 true cm

\vskip 0.3 true cm


\section{Introduction}
The 3D incompressible Navier-Stokes(NS) equations are given as
\be\label{NS}
\lt\{
\begin{aligned}
&\p_t u+u\cdot\na u+\na p-\Dl u=0,\\
&\na\cdot u=0,
\end{aligned}
\rt.
\ee
where $u(x,t)\in\bR^3,p(x,t)\in\bR$ represents the velocity vector and the scalar pressure. The NS equations are one of the most fundamental nonlinear partial differential equations in nature but are far from being fully understood. The global regularity problem of solutions to the 3D NS equations with smooth initial data remains open and is viewed as one of the most important open questions in mathematics \cite{Fcl1}.

To study the possible singularity formation for the solution of NS equations, one often scales the solution at the possible singular point. This results in a nontrivial bounded solution existing in the whole space $\bR^3$ and the time interval $(-\i,0]$, which is often referred to as an ancient solution. Information of the ancient solution reveals the singular structure of the original solution. In some sense, the trivialness($u\equiv 0$) of the ancient solution equals to the regularity of the solution of NS equations. To study the Liouville type theorem of the ancient solution to the NS equations seems to be the first step to understand the regularity of the solution of NS equations. However, the Liouville theorem of ancient solutions seems also beyond touch if no extra assumption is given. In fact, it is still wildly open even for the stationary case since another old unsolved problem concerning D-solution, which asks if a 3D steady solution of NS equations with finite Dirichlet integral and vanishing at infinity is zero identically. This problem is not solved even in the axially symmetric case. See some recent paper, for examples, \cite{CPZZ:2, Sg:1}.

In 1934, Jean Leray \cite{Lj:1} raised the existence of back self-similar solutions of NS equations, which can be viewed as ancient solutions with a uniform profile. In \cite{NRS:1}, the authors proved that such solutions must be trivial if the profile belongs to $L^p(p=3)$. Later Tsai \cite{Tsai:1998ARMA} improves this result to the case $p\in (3,+\i]$. See also an extension in Chae \cite{CJ:2017ARMA}. In the remarkable paper \cite{ISS:2003RMS}, the authors proved that $L^\i L^3$ solution of NS equations must be regular. The above four papers are all based on the landmark partial regularity theory of Caffarelli-Kohn-Nirenberg \cite{CKN:1982CPAM}. In \cite{CKN:1982CPAM}, the authors showed that the 1-dimensional Hausdorff measure of the singular set of suitable weak solutions must be zero, which indicates that for the axially symmetric Navier-Stokes(ASNS) equations, blow up can only happen on the symmetric axis.

In this paper, we consider the ASNS equations. In the cylindrical coordinates $(r,\th, z)$, we have $x=(x_1,x_2,x_3)=(r\cos\th,r\sin\th,z)$ and the axially symmetric solution of the incompressible Navier-Stokes equations is given as
\[
u=u^r(r,z,t)e_r+u^{\th}(r,z,t)e_{\th}+u^z(r,z,t)e_z,
\]
where the basis vectors $e_r,e_\th,e_z$ are
\[
e_r=(\frac{x_1}{r},\frac{x_2}{r},0),\quad e_\th=(-\frac{x_2}{r},\frac{x_1}{r},0),\quad e_z=(0,0,1).
\]

The components $u^r,u^\th,u^z$, which are independent of $\th$, satisfy
\begin{equation}\label{ASNS}
\left\{
\begin{aligned}
&\p_t u^r+(u^r\p_r+u^z\p_z)u^r -\frac{(u^\th)^2}{r}+\p_r p=\left(\Delta-\frac{1}{r^2}\right)u^r, \\
&\p_tu^\th+(u^r\p_r+u^z\p_z) u^\th+\frac{u^\th u^r}{r}=\left(\Delta-\frac{1}{r^2}\right)u^\th , \\
&\p_t u^z+(u^r\p_r+u^z\p_z)u^z+\p_z p=\Delta u^z ,                                    \\
&\p_ru^r+\frac{u^r}{r}+\p_zu^z=0.
\end{aligned}
\right.
\end{equation}

We can also compute the axi-symmetric vorticity $w=\nabla\times u=w^re_r+w^\th e_\th+w^ze_z$  as follows
\[
w^r=-\p_z u^\th, \ w^\th=\p_z u^r-\p_r u^z,\  w^z=\left(\p_r+\frac{1}{r}\right)u^\th,
\]
which satisfies
\be\label{VEQ}
\lt\{
\begin{aligned}
&\p_t w^r+(u^r\p_r+u^z\p_z)w^r-\left(\Dl-\f{1}{r^2}\right)w^r-(w^r\p_r+w^z\p_z)u^r=0,\\
&\p_t w^\th+(u^r\p_r+u^z\p_z)w^\th-\left(\Dl-\f{1}{r^2}\right)w^\th-\f{u^r}{r}w^\th-\f{1}{r}\p_z(u^\th)^2=0,\\
&\p_t w^z+(u^r\p_r+u^z\p_z)w^z-\Dl w^z-(w^r\p_r+w^z\p_z)u^z=0.
\end{aligned}
\rt.
\ee
Global in-time regularity of the solution to the ASNS equations is still open. In \cite{CSTY:2}, Chen-Strain-Yau-Tsai proved that a suitable weak solution is regular if the solution satisfies $r|b|\leq C_\ast<\i$. Also, Koch-Nadirashvili-Seregin-Sverak in \cite{KNSS:1} proved the same result by using a Liouville theorem and the scaling-invariant property of NS equations. Lei-Zhang in \cite{LZ:2011JFA} proved regularity of the solution under a more general assumption on the drift term $b$ where $b\in L^\i\left([-1,0),BMO^{-1}\right)$.

In the paper \cite{KNSS:1}, the authors make the following conjecture: {\it any bounded ancient mild solution of the 3D axially symmetric Navier-Stokes equations is constant.} And it is proved under the assumption that the solution is swirl free ($u^\th\equiv 0$) in \cite{KNSS:1}. Our purpose of this paper is to improve their result by allowing that the solution can grow with a power smaller than 1 with respect to the distance to the origin. Also, we will show that such a power is optimal to prove the Liouville type theorem since we can find counterexamples for the Navier-Stokes equations such that the Liouville theorem fails if the solution can grow linearly.

Under the situation that $u^\th\equiv 0$, the ancient solution of ASNS reads
\begin{equation}\label{ASNS1}
\left\{
\begin{aligned}
&\p_t u^r+(u^r\p_r+u^z\p_z)u^r+\p_r p=\left(\p^2_r+\f{1}{r}\p_r+\p^2_z-\frac{1}{r^2}\right)u^r,\q \text{in}\ (-\i,0]\times\bR^3,\\
&\p_t u^z+(u^r\p_r+u^z\p_z)u^z+\p_z p=\left(\p^2_r+\f{1}{r}\p_r+\p^2_z\right) u^z ,              \qq\q\ \text{in}\ (-\i,0]\times\bR^3,         \\
&\p_ru^r+\frac{u^r}{r}+\p_zu^z=0,                                          \qq\qq\qq\q\qq\qq\qq\qq\ \ \ \ \ \text{in}\ (-\i,0]\times\bR^3.
\end{aligned}
\right.
\end{equation}
Set $\O=\f{w^\th}{r}$, which satisfies
\be\label{eo}
\Delta\Omega+\frac{2}{r}\p_r\Omega-\p_t\Omega=(u^r\p_r+u^z\p_z)\Omega.
\ee
Here is the main result of this paper:
\begin{theorem}\label{th1}
Suppose the ancient solution $u$ of ASNS equations is smooth and no-swirl, then $u=(u^r,u^z)=(0,c(t))$ provided that
\be\label{egc}
|u^r(t,r,z)|+|u^z(t,r,z)|\leq C(\s{-t}+|x|)^\al,
\ee
\be\label{egc1}
|u^r(t,r,z)|+|u^z(t,r,z)-u^z(t,0,z)|\leq Cr,
\ee
for any $\alpha<1$. Here $C$ is a constant independent of $t$, $r$ and $z$.
\end{theorem}

Some remarks are in followings.
\begin{remark}
The assumption \eqref{egc} indicates that $u$ is sublinearly growing at infinity, which corresponds to the maximally allowed growing condition for the Liouville theorem for the heat equation.
\end{remark}
\begin{remark}
Since the axially symmetic solution is smooth, from \cite{LW:2009JMA}, we have
\be
\p^{2l}_r u^r|_{r=0}=\p^{2l+1}_r u^z|_{r=0}=0,\q \text{for} \q l\in \{0\}\cup\bN.\nn
\ee
So, when $r<1$, we can write
\be
u^r(t,r,z)=rf(t,r,z),\q u^z(t,r,z)-u^z(t,0,z)=rg(t,r,z),\nn
\ee
where $f,g$ are smooth for $t,r,z$. Our assumption \eqref{egc1} actually is the restriction that $f,g$ are uniformly bounded with respect to $t,z$.
\end{remark}
\begin{remark}
For 2-dimensional steady Navier-Stokes equations, the authors in \cite{BFZ:2013JMFM} proved a Liouville-type theorem by assuming that $\limsup_{|x|\rightarrow \i}|u(x)||x|^{-\al}< \i$ for some $\al<1/3$. For 2-dimensional ancient mild solution of Navier-Stokes equations, the authors in \cite{LZZ:2017} derived a Liouville-type theorem under the condition that $\lim_{|x|\rightarrow \i}\left(|u(t,x)||x|^{-1}+|w(t,x)|\right)=0$. They also proved a Liouville-type theorem of 3-dimensional ASNS system without swirl under the conditions $\lim_{r\rightarrow \i}\O=0$ and $\lim_{|x|\rightarrow \i}|u(t,x)||x|^{-1}=0$. The main difference between the result in \cite{LZZ:2017} and ours is that we do not impose any assumption on the vorticity.
\end{remark}

Now if $\al=1$ in the assumption \eqref{egc}, our Liouville-type theorem \eqref{th1} fails. Consider the ASNS equations without swirl. We have the following proposition:
\begin{proposition}\label{pro}
If $u$ is the solution of ASNS equation and satisfies
\be
|u^r(t,r,z)|+|u^z(t,r,z)-u^z(t,0,z)|\leq C_\ast r.\nn
\ee
Then we have the following a family of solutions
\be\label{ess}
\begin{array}{c}
u^r=C_1r,\ u^z=-2C_1z+C_2(t),\\
p=-C^2_1\left(\f{1}{2}r^2+2z^2\right)+(2C_1C_2(t)-C'_2(t))z+C_3,
\end{array}
\ee
for any constants $C_1\leq C_\ast$, $C_3$, and smooth function $C_2(t)$.
\end{proposition}
\begin{remark}
When $C_1\neq 0$, the solution \eqref{ess} with $C_2(t)=0$ satisfies the assumption \eqref{egc} with $\al=1$, which indicates the assumption \eqref{egc} for the power $\al<1$ is sharp to prove the Liouville-type theorem. We also note that the kind of linear solutions have also been observed in \cite{LZZ:2017}, see Remark 4.1.
\end{remark}
Now we outline the proof of Theorem \ref{th1} briefly. We will use \eqref{eo} to show that for sufficiently large $q$, we have $\|\O\|_{L^q((-\i,0]\times\bR^3)}=0$, which indicates that $w^\th\equiv 0$. Then by using $Biot-Savart$ law, we can show that the solution satisfies the Laplace's equation which implies that $u$ is a vector depending only on time if it is sublinearly growing with distance to the origin.

Throughout this paper, $C$ denotes a positive constant which may be different from line to line. Meanwhile, we denote
\be
D_R=\left\{(r,z):\,0\leq r\leq R,\,|z|\leq R\right\}.
\ee
We also apply $A\lesssim B$ to denote $A\leq CB$.
\section{Proof of Theorem \ref{th1}}

\noindent{\bf Vanishing of $\boldsymbol{w^\th}$}\\
Our main procedure is to show that, for sufficiently large $q\geq 2$,
\be\label{eo1}
\|\O\|_{L^q((-\i,0]\times\bR^3)}=0,
\ee
with the help of equation \eqref{eo} and assumptions \eqref{egc}, \eqref{egc1}.
Let $q\in\mathbb{R}_+$ be determined later and $\eta$ be a cut-off function such that
\be
\eta=\eta(t,r,z)=\eta_0(t)\eta_1(r)\eta_2(z)
\ee
with
\be
\begin{footnotesize}
\eta_0(t)=\left\{
\begin{aligned}
&1,\quad\text{if }\quad t\in[-1,0]; \\
&0,\quad\text{if }\quad t\leq -2, \\
\end{aligned}
\right.
\ \eta_1(r)=\left\{
\begin{aligned}
&1,\quad\text{if }\quad r\leq 1; \\
&0,\quad\text{if }\quad r\geq 2, \\
\end{aligned}
\right.
\ \eta_2(z)=\left\{
\begin{aligned}
&1,\quad\text{if }\quad |z|\leq 1; \\
&0,\quad\text{if }\quad |z|\geq 2, \\
\end{aligned}
\right.
\end{footnotesize}
\ee
and $|\eta'_0,\,\eta_1',\,\eta_2'|\leq2$. Meanwhile, $\eta_R$ is denoted by
\be
\eta_R(t,r,z):=\eta\left(\frac{t}{R^2},\,\frac{r}{R},\frac{z}{R}\right),
\ee
and it is easy to see that
\be\label{etestf}
|\p_t\eta_R|+|\na\eta_R|^2\leq CR^2.
\ee
Direct calculation shows
\bea\label{3.6}
&&\q\|\O\eta_R\|^q_{L^q((-\i,0]\times\bR^3)}\nn\\
&&=2\pi\int_{-\infty}^0\int_{-\infty}^\infty\int_0^\infty|\Omega\eta_R|^{q}rdrdzdt\nn\\
&&=2\pi\int_{-\infty}^0\int_{-\infty}^\infty\int_0^\infty\left(\p_zu^r-\p_ru^z\right)\Omega|\Omega|^{q-2}\eta_R^{q}drdzdt\nn\\
&&=2\pi\int_{-\infty}^0\int_{-\infty}^\infty\int_0^\infty\Big\{\p_zu^r-\p_r\Big(u^z-u^z\Big|_{r=0}\Big)\Big\}\Omega|\Omega|^{q-2}\eta^{q}_Rdrdzdt\nn\\
&&=-2\pi\int_{-\infty}^0\int_{-\infty}^\infty\int_0^\infty u^r\p_z\left(\Omega|\Omega|^{q-2}\right)\eta^{q}_Rdrdzdt\nn\\
&&\hskip 1cm+2\pi\int_{-\infty}^0\int_{-\infty}^\infty\int_0^\infty \Big(u^z-u^z\Big|_{r=0}\Big)\p_r\left(\Omega|\Omega|^{q-2}\right)\eta^{q}_Rdrdzdt\nn\\
&&\hskip 1cm-2\pi\int_{-\infty}^0\int_{-\infty}^\infty\int_0^\infty u^r\Omega|\Omega|^{q-2}\p_z\eta_R^{q}drdzdt\nn\\
&&\hskip 1cm+2\pi\int_{-\infty}^0\int_{-\infty}^\infty\int_0^\infty \Big(u^z-u^z\Big|_{r=0}\Big)\Omega|\Omega|^{q-2}\p_r\eta_R^{q}drdzdt\nn\\
&&\leq (q-1)\int_{-\infty}^0\int_{\bR^3}\left|\frac{u^r}{r}\right|\cdot|\p_z\Omega|\cdot|\Omega|^{q-2}\eta^{q}_Rdxdt\nn\\
&&\hskip 1cm +(q-1)\int_{-\infty}^0\int_{\bR^3}\left|\frac{u^z-u^z\big|_{r=0}}{r}\right||\p_r\Omega|\cdot|\Omega|^{q-2}\eta^{q}_Rdxdt\nn\\
&&\hskip 1cm+q\int_{-\infty}^0\int_{\bR^3}\left|\frac{u^r}{r}\right|\cdot|\Omega|^{q-1}\eta^{q-1}_R|\p_z\eta_R|dxdt\nn\\
&&\hskip 1cm +q\int_{-\infty}^0\int_{\bR^3}\left|\frac{u^z-u^z\big|_{r=0}}{r}\right||\Omega|^{q-1}\eta^{q-1}_R|\p_r\eta_R|dxdt,
\eea
since $u^r(t,0,z)\equiv0$. By the condition \eqref{egc}, it follows that
\be
\left|\f{u^r}{r}\right|+\left|\f{u^z-u^z\big|_{r=0}}{r}\right|\leq C.\nn
\ee
Then we have
\be\label{3.7}
\begin{split}
&\q\|\O\eta_R\|^q_{L^q((-\i,0]\times\bR^3)}\\
&\leq C_q\int_{-\infty}^0\int_{\bR^3}|\bar{\na}\Omega||\Omega|^{q-2}\eta^{q}_Rdxdt+C_q\int_{-\infty}^0\int_{\bR^3}|\Omega|^{q-1}\eta^{q-1}_R|\bar{\na}\eta_R|dxdt,\\
\end{split}
\ee
where $\bar{\na}=(\p_r,\p_z)$.
Using Young inequality, \eqref{3.7} follows that
\be\label{2.6}
\begin{split}
&\q\|\O\eta_R\|^q_{L^q((-\i,0]\times\bR^3)}\\
&\leq \f{1}{4}\int_{-\infty}^0\int_{\bR^3}|\Omega|^{q}\eta^{q}_Rdxdt+C_q\int_{-\infty}^0\int_{\bR^3}|\bar{\na}\Omega|^2|\Omega|^{q-4}\eta^{q}_Rdxdt\\
&\q+\f{1}{4}\int_{-\infty}^0\int_{\bR^3}|\Omega|^{q}\eta^{q}_Rdxdt+C_q\int_{-\infty}^0\int_{\bR^3}|\bar{\na}\eta_R|^qdxdt.\\
\end{split}
\ee

The first and third terms on the right hand of \eqref{2.6} can be absorbed by the left hand side. Then one derive
\bea\label{2.8}
&&\q\|\O\eta_R\|^q_{L^q((-\i,0]\times\bR^3)}\nn\\
&\leq& C_q\left(\int_{-\infty}^0\int_{\bR^3}|\bar{\na}\Omega|^2|\Omega|^{q-4}\eta^{q}_Rdxdt+\int_{-\infty}^0\int_{\bR^3}|\bar{\na}\eta_R|^qdxdt\right)\nn\\
&\leq& C_q\int_{-\infty}^0\int_{\bR^3}|\bar{\na}\Omega|^2|\Omega|^{q-4}\eta^{q}_Rdxdt+C_qR^{-q}\int_{\supp\eta_R}dxdt\nn\\
&\leq& C_q\underbrace{\int_{-\infty}^0\int_{\bR^3}|\bar{\na}\Omega|^2|\Omega|^{q-4}\eta^{q}_Rdxdt}_I+C_qR^{-q+5},
\eea
where we have applied \eqref{etestf}.

Now we come to estimate $I$ on the right hand of \eqref{2.8}. Testing \eqref{eo} by $\Omega|\Omega|^{q-4}\eta_R^{q}$ indicates that
\bea\label{2.9}
&&\int^0_{-\i}\int_{\bR^3}\Big(\p_t-\Delta-\f{2}{r}\p_r\Big)\O\cdot\Omega|\Omega|^{q-4}\eta_R^{q}dxdt\nn\\
&=&-\int^0_{-\i}\int_{\bR^3}(u^r\p_r+u^z\p_z)\Omega \cdot\Omega|\Omega|^{q-4}\eta_R^{q}dxdt.
\eea
Integration by parts implies
\be\label{2.10}
\begin{split}
&\f{1}{q-2}\int_{\bR^3}|\O(0,x)|^{q-2}\eta_R^q(0,x)dx+(q-3)\int_{-\infty}^0\int_{\bR^3}|\nabla\Omega|^2|\Omega|^{q-4}\eta_R^{q}dxdt\\
&+\frac{4\pi}{q-2}\int_{-\infty}^0\int_{-\infty}^\infty|\Omega(t,0,z)|^{q-2}\eta^{q}_R(t,0,z)dzdt\\
=&-\int_{-\infty}^0\int_{\bR^3}\nabla\Omega\cdot\nabla\eta_R^{q}\Omega|\Omega|^{q-4}dxdt-\frac{4\pi}{q-2}\int_{-\infty}^0\int_{-\infty}^\infty\int_0^\infty|\Omega|^{q-2}\p_r\eta^{q}_Rdrdzdt\\
&+\frac{1}{q-2}\int_{-\infty}^0\int_{\bR^3} u\cdot\nabla\eta^{q}_R|\Omega|^{q-2}dxdt+\frac{1}{q-2}\int_{-\infty}^0\int_{\bR^3} |\Omega|^{q-2}\p_t\eta^{q}_R dxdt.
\end{split}
\ee
By Young inequality, the first term on the right hand side of \eqref{2.10} could be estimated by
\be\label{3.9}
\begin{split}
&\left|\int_{-\infty}^0\int_{\bR^3}\nabla\Omega\cdot\nabla\eta_R^{q}\Omega|\Omega|^{q-4}dxdt\right|\\
\leq&\,\frac{q-3}{2}\int_{-\infty}^0\int_{\bR^3}|\nabla\Omega|^2\cdot|\Omega|^{q-4}\eta_R^{q}dxdt\\
&+C_q\int_{-\infty}^0\int_{\bR^3}|\Omega|^{q-2}|\nabla\eta_R|^2\eta_R^{q-2}dxdt.
\end{split}
\ee
Combining \eqref{2.10} and \eqref{3.9}, it follows that
\be\label{2.12}
\begin{split}
I=&\int_{-\infty}^0\int_{\bR^3}|\nabla\Omega|^2\cdot|\Omega|^{q-4}\eta_R^{q}dxdt\\
\leq&\,C_q\left(\int_{-\infty}^0\int_{\bR^3}|\Omega|^{q-2}\left(|\nabla\eta_R|^2\eta_R^{q-2}+|\p_r\eta_R|\eta_R^{q-1}r^{-1}+|\p_t\eta_R|\eta_R^{q-1}\right)dxdt\right.\\
&\left.+\|u\|_{L^\infty([-4R^2,0]\times D_{2R})}\int_{-\infty}^0\int_{\bR^3}|\Omega|^{q-2}\cdot|\nabla\eta_R|\eta_R^{q-1}dxdt\right)\\
=&J_1+J_2+J_3+J_3.
\end{split}
\ee
Using Young inequality and \eqref{etestf}, we have
\bea\label{2.13}
J_1&\leq& \ve^{\f{q}{q-2}}\int_{-\infty}^0\int_{\bR^3}|\Omega\eta_R|^qdxdt+\ve^{-\f{q}{2}}C_q\int_{-\infty}^0\int_{\bR^3}|\na\eta_R|^qdxdt\nn\\
&\leq& \ve^{\f{q}{q-2}}\int_{-\infty}^0\int_{\bR^3}|\Omega\eta_R|^qdxdt+\ve^{-\f{q}{2}}C_qR^{-q+5}.
\eea
And
\bea\label{2.14}
J_2&\leq& \ve^{\f{q}{q-2}}\int_{-\infty}^0\int_{\bR^3}|\Omega\eta_R|^qdxdt+\ve^{-\f{q}{2}}C_q\int_{-\infty}^0\int_{\bR^3}|\p_r\eta_R\eta_Rr^{-1}|^{q/2}dxdt\nn\\
&\leq& \ve^{\f{q}{q-2}}\int_{-\infty}^0\int_{\bR^3}|\Omega\eta_R|^qdxdt+\ve^{-\f{q}{2}}C_qR^{-q+5}.
\eea
Also
\bea\label{2.15}
J_3&\leq& \ve^{\f{q}{q-2}}\int_{-\infty}^0\int_{\bR^3}|\Omega\eta_R|^qdxdt+\ve^{-\f{q}{2}}C_q\int_{-\infty}^0\int_{\bR^3}|\p_t\eta_R\eta_R|^{\f{q}{2}}dxdt\nn\\
&\leq& \ve^{\f{q}{q-2}}\int_{-\infty}^0\int_{\bR^3}|\Omega\eta_R|^qdxdt+\ve^{-\f{q}{2}}C_qR^{-q+5}.
\eea
At last, By using \eqref{egc}, we have
\be
\|b\|_{L^\infty([-4R^2,0]\times D_{2R})}\lesssim R^{\alpha},\quad\text{with }\alpha<1.
\ee
Then we have
\bea\label{2.17}
J_4&\leq& \ve^{\f{q}{q-2}}\int_{-\infty}^0\int_{\bR^3}|\Omega\eta_R|^qdxdt+\ve^{-\f{q}{2}}C_qR^{\f{q}{2}\al}\int_{-\infty}^0\int_{\bR^3}|\na\eta_R\eta_R|^{q/2}dxdt\nn\\
&\leq& \ve^{\f{q}{q-2}}\int_{-\infty}^0\int_{\bR^3}|\Omega\eta_R|^qdxdt+\ve^{-\f{q}{2}}C_qR^{\f{q}{2}(\al-1)+5}.
\eea
Combining the estimates from \eqref{2.12} to \eqref{2.17}, we can obtain
\be\label{2.18}
I\leq  4\ve^{\f{q}{q-2}} \int_{-\infty}^0\int_{\bR^3}|\Omega\eta_R|^qdxdt+\ve^{-\f{q}{2}}C_q\big(R^{\f{q}{2}(\al-1)+5}+R^{-q+5}\big).
\ee
Remembering \eqref{2.8}, substituting \eqref{2.18} into it and by choosing sufficiently small $\ve$, we can obtain
\bea
\q\|\O\eta_R\|^q_{L^q((-\i,0]\times\bR^3)}\leq C_q\big(R^{\f{q}{2}(\al-1)+5}+R^{-q+5}\big).
\eea
By choosing $q>\frac{10}{1-\alpha}$ and then let $R\to+\infty$, it follows that
\be
\int_{-\infty}^0\int_{\bR^3}|\Omega|^{q}dxdt=0,\quad\forall q>\frac{10}{1-\alpha}.
\ee
Therefore we conclude that $\Omega\equiv0$. Then we can get $w^\th\equiv 0$.

\noindent{\bf Liouville theorem of $\boldsymbol{u}$}\\
For the axially symmetric solution without swirl, we have
\be
\na\times (u^re_r+u^ze_z)=w^\th e_\th=0.\nn
\ee From $Biot-Savart$ law,
\be
\na\times\na\times u=-\Dl u+\na(\na\cdot u),\nn
\ee
and divergence-free condition of $u$, we have $\Dl u=0$.

Since $u$ satisfies the assumption \eqref{egc}, which is sublinearly growing with respect to $|x|$, we can get
\be
u=(u^r,u^z)=(c_1(t),c_2(t)).\nn
\ee
Also we have $|u^r|\leq Cr$, so $c_1(t)=0$. This finishes the proof of Theorem \ref{th1}.
\ef
\section{Proof of the Proposition \ref{pro}}
This is a direct computation, so we omit the details.\ef

\section*{Acknowledgments} The authors wish to thank Prof. Qi S. Zhang in UC Riverside for his constant encouragement and courtesy reminding of related results on this topic. X. Pan is supported by Natural Science Foundation of Jiangsu Province (No. SBK2018041027) and National Natural Science Foundation of China (No. 11801268). Z. Li is supported by the Startup Foundation for Introducing Talent of NUIST (No. 2019r033).

\end{document}